\documentclass[12pt]{amsart}
\usepackage{amssymb}
\usepackage{amsmath}
\usepackage{graphicx}
\newtheorem{theorem}{Theorem}[section]

\newtheorem{lemma}{Lemma}[section] 
\newtheorem{corollary}{Corollary}[section]
\newtheorem{exa}[theorem]{Example}

\theoremstyle{definition}

\newtheorem{definition}{Definition}[section]

\begin{document}
\title{\large Combinations and Mixtures of Optimal Policies in Unichain Markov Decision Processes are Optimal}
\author{\bf Ronald ORTNER}\thanks{~This work was supported in 
part by the the Austrian Science Fund FWF (S9104-N04 SP4) and the IST Programme of the European 
Community, under the PASCAL Network of Excellence, IST-2002-506778. This publication only reflects 
the authors' views.}
\email{rortner@unileoben.ac.at}
\address{\parbox{1.4\linewidth}{Department Mathematik und Informationstechnolgie\\ Montanuniversit\"at Leoben\\
Franz-Josef-Strasse 18\\8700 Leoben, Austria}}
\begin{abstract}
We show that combinations of optimal (stationary) policies in unichain Markov decision processes are optimal. 
That is, let  $\mathcal M$ be a unichain Markov decision process with state space $S$, action space $A$ and 
policies $\pi_j^\circ: S\to A$ ($1\leq j\leq n$) with optimal  average infinite horizon reward. Then any  
\textit{combination} $\pi$ of these policies, where for each state $i\in S$ there is a $j$ such that 
$\pi(i)=\pi_j^\circ(i)$, is optimal as well. Furthermore, we prove that any \textit{mixture} of optimal policies,
where at each visit in a state $i$ an arbitrary action $\pi_j^\circ(i)$ of an optimal policy is chosen, yields 
optimal average reward, too.
\end{abstract} 
\maketitle

\markboth{\sf R. Ortner}{\sf Combinations and Mixtures of Optimal Policies in Unichain MDPs are Optimal}
\baselineskip 15pt

\section{Introduction}
\begin{definition}
A \emph{Markov decision process} (MDP) $\mathcal{M}$ on a (finite) set of \emph{states} $S$ with a 
(finite) set of \emph{actions} $A$ available in each state $\in S$ consists of 
\begin{enumerate}
\renewcommand{\theenumi}{(\roman{enumi})}
\item an initial distribution $\mu_0$ that specifies the probability of starting in some state in $S$,
\item the transition probabilities $p_a(i,j)$ that specify 
the probability of reaching state $j$ when choosing action $a$ in state $i$, and 
\item the payoff distributions with mean $r_a(i)$ that 
specify the random reward for choosing action $a$ in state $i$.
\end{enumerate}

A (stationary) \emph{policy} on $\mathcal M$ is a mapping $\pi: S \to A$.
\end{definition}

Note that each policy $\pi$ induces a Markov chain on $\mathcal{M}$. We are interested in MDPs, where 
in each of the induced Markov chains any state is reachable from any other state.

\begin{definition}
An MDP $\mathcal M$ is called \textit{unichain}, if for each policy $\pi$ the Markov
chain induced by $\pi$ is ergodic, i.e. if the matrix $P=(p_{\pi(i)}(i,j))_{i,j\in S}$ is irreducible.
\end{definition}

It is a well-known fact (cf.\ e.g.\ \cite{keme}, p.130ff) that for an ergodic Markov chain with 
transition matrix $P$ there exists a unique invariant and strictly positive distribution $\mu$, such that 
independent of the initial distribution $\mu_0$ one has $\mu_n=\bar{P}_n\mu_0\to \mu$, 
where $\bar{P}_n=\frac{1}{n}\sum_{j=1}^nP^j$.\footnote{Actually, for aperiodic Markov chains 
one has even $P^n\mu_0\to\mu$, while the convergence behavior of periodic Markov chains 
can be described more precisely. However, for our purposes the stated fact is sufficient.}
Thus, given a policy $\pi$ on a unichain MDP that induces a Markov chain with invariant 
distribution $\mu$, the \textit{average reward} of that policy can be defined as
$$V({\pi}):=\sum_{i\in S} \mu(i) r_{\pi(i)}(i).$$
A policy $\pi^\circ$ is called \textit{optimal} if for all policies $\pi$: $V({\pi})\leq V({\pi^\circ}).$ It can be shown 
(\cite{pute}, p.360ff) that in the unichain case the optimal value $V({\pi^\circ})$ cannot be increased by allowing 
time-dependent policies, as there is always a \textit{stationary} (time-independent) policy that 
gains optimal average reward, which is why we consider only stationary policies.

In this setting we are going to prove that \textit{combinations} of optimal policies are optimal as well.
\begin{theorem}\label{theorem}
Let $\mathcal M$ be a unichain MDP with state space $S$ and $\pi_1^{\circ},\pi_2^{\circ}$ optimal 
policies on $\mathcal M$. Then any  \textit{combination} $\pi$ of these policies where for each state 
$i\in S$ either $\pi(i)=\pi_1^\circ (i)$ or $\pi(i)=\pi_2^\circ (i)$ is optimal as well.
\end{theorem}

Obviously, if two combined optimal policies are optimal, so are combinations of an arbitrary number of optimal
policies. Thus, one immediately obtains that the set of optimal policies is closed under combination.

\begin{corollary}
Let $\mathcal M$ be a unichain MDP with state space $S$. A policy $\pi$ is optimal on $\mathcal M$ if and only if
for each state $i$ there is an optimal policy $\pi^\circ$ with $\pi(i)=\pi^\circ (i)$.
\end{corollary}

\section{Proof of Theorem \ref{theorem}}
We start with a result about the distributions of policies that differ in at most two states.
\begin{lemma}\label{lem:distr}
Let $\mathcal M$ be a unichain MDP with state space $S$. Let $\pi_{00},\pi_{01},\pi_{10},\pi_{11}$ be four policies on   $\mathcal M$ with invariant distributions $\mu_{00}=(a_i)_{i\in S}$, $\mu_{01}=(b_i)_{i\in S}$, 
$\mu_{10}=(c_i)_{i\in S}$, $\mu_{11}=(d_i)_{i\in S}$ and $s_1,s_2$ two states in $S$ such that 
\begin{enumerate}
\renewcommand{\theenumi}{(\roman{enumi})}
\item for all $i\in S\setminus\{s_1,s_2\}$: $\pi_{00}(i)=\pi_{01}(i)=\pi_{10}(i)=\pi_{11}(i)$,
\item  $\pi_{00}(s_1)=\pi_{01}(s_1)\neq \pi_{10}(s_1)=\pi_{11}(s_1)$,
\item $\pi_{00}(s_2)=\pi_{10}(s_2)\neq \pi_{01}(s_2)=\pi_{11}(s_2)$.
\end{enumerate}
Then each of the distributions $\mu_{ij}$ is uniquely determined by the other three. More precisely, 
e.g.\footnote{For the sake of readability we don't give a general formula but only reproduce how to
calculate $\mu_{11}$. Since the situation is symmetric it is easy to see (but a bit tedious to write down
or read) how the general formula looks like.} for all states $i$
\[ d_i=\frac{a_{s_2}b_{s_1}c_i - a_i b_{s_1} c_{s_2} + a_{s_1}b_i c_{s_2}}
		{a_{s_2}b_{s_1}- b_{s_1} c_{s_2} +  a_{s_1} c_{s_2}}.
\]
\end{lemma}
\begin{proof}
Since $\mathcal M$ is unichain, the distributions $\mu_{ij}$ are all uniquely determined by the transition
matrices $P_{ij}$ of the Markov chains induced by the policies $\pi_{ij}$. By assumption (i), the matrices $P_{ij}$
share all rows except rows $s_1,s_2$, which we may assume to be the first and second row, respectively. 
Furthermore, by (ii), $P_{00}$ and $P_{01}$ share the first row as well as $P_{10}$ and $P_{11}$. Finally, 
by (iii) we have equal second rows in $P_{00}$ and $P_{10}$ as well as in $P_{01}$ and $P_{11}$. 

Since by assumption the distributions are invariant, we have $\mu_{ij}P_{ij}=\mu_{ij}$. 
Writing the probabilities in $P_{00}$ as $p_{ij}$ and those in the first two rows of $P_{11}$
as $q_{ij}$,  it follows that for each state $i$:
\begin{eqnarray}   
a_i=\sum_{j\in S} a_jp_{ji}, \qquad  b_i= b_2q_{2i} + \!\!\sum_{j\in S\setminus\{2\}} b_jp_{ji},\qquad 
c_i=c_1q_{1i} + \!\!\sum_{j\in S\setminus\{1\}} c_jp_{ji}. \label{eq1}
\end{eqnarray}
Setting $\nu=(\nu_i)_{i\in S}$ with $\nu_i:=a_2b_1c_i-a_ib_1c_2+a_1b_ic_2$, one has by (\ref{eq1})
\begin{eqnarray*} 
(\nu P_{11})_i &=& a_2b_1c_1q_{1i} + a_1b_2c_2q_{2i}
+ \sum_{j\in S\setminus\{1,2\}} a_2b_1c_jp_{ji} \\
&&\qquad\qquad-\sum_{j\in S\setminus\{1,2\}} a_jb_1c_2p_{ji} 
+ \sum_{j\in S\setminus\{1,2\}} a_1b_jc_2p_{ji}\\
&=& a_2b_1c_1q_{1i} + a_1b_2c_2q_{2i}
+ a_2b_1 \Big(c_i-c_1q_{1i}-c_2p_{2i}\Big)\\
&&- b_1c_2  \Big(a_i-a_1p_{1i}-a_2p_{2i}\Big)
+ a_1c_2 \Big(b_i-b_1p_{1i}-b_2q_{2i}\Big)\\
&=&a_2b_1c_i-a_ib_1c_2+a_1b_ic_2=\nu_i
\end{eqnarray*}

Hence, normalizing $\nu$ one has an invariant distribution of $P_{11}$, which by assumption is unique
and consequently identical to $\mu_{11}$.
\end{proof}

With this information on the distributions, we are able to tell something about the average rewards of
the policies as well.

\begin{lemma}\label{lem:comp}
Let $\pi_{00},\pi_{01},\pi_{10},\pi_{11}$ and $s_1,s_2$ be as in Lemma \ref{lem:distr} and denote the average
rewards of the policies by $V_{00},V_{01},V_{10},V_{11}$. Let $a,b\in\{0,1\}$ and set $\neg x:=1-x$. Then 
it cannot be the case that $V_{ab}>V_{\neg a, b},V_{a, \neg b}$ and $V_{\neg a, \neg b}\geq V_{\neg a, b},V_{a, \neg b}$.
Analogously, it cannot hold that $V_{ab}<V_{\neg a, b},V_{a, \neg b}$ and 
$V_{\neg a, \neg b}\leq V_{\neg a, b},V_{a, \neg b}$.
\end{lemma}

\begin{proof}
For the sake of readability, we prove the case $a=b=0$. The other cases follow by symmetry. Actually, we will show 
that if $V_{00}>V_{01},V_{10}$, then $V_{00}>V_{11}$. Since one has analogously the implication
that if $V_{11}\geq V_{01},V_{10}$, then  $V_{11}\geq V_{00}$, the assumptions $V_{00}>V_{01},V_{10}$
and $V_{11}\geq V_{01},V_{10}$ obviously lead to a contradiction.

Similarly as in the case of transition probabilities, in the following we write for the rewards of the 
policy $\pi_{00}$ simply $r_i$ instead of $r_{\pi_{ij}(i)}(i)$. For the deviating rewards in state 
$s_1$ under policies $\pi_{10},\pi_{11}$ and state $s_2$ under $\pi_{01},\pi_{11}$ we write
$r'_1$ and $r'_2$, respectively. Then we have
\begin{align*}  
&V_{00}=\sum_{i\in S} a_ir_i,
&V_{01}=& b_2r'_2 + \!\!\sum_{i\in S\setminus\{2\}} b_ir_i,\\
&V_{10}=c_1r'_1 + \!\!\sum_{i\in S\setminus\{1\}}\!\! c_ir_i, 
&V_{11}=&d_1r'_1 +d_2r'_2 +\!\!\sum_{i\in S\setminus\{1,2\}}\!\! d_ir_i. 
\end{align*}
If we now assume that $V_{00}>V_{01},V_{10}$, the first three equations yield
\begin{eqnarray}  
b_2r'_2&<& a_2r_2 + \sum_{i\in S\setminus\{2\}} (a_i-b_i)r_i,\label{eq2}\\
c_1r'_1&<& a_1r_1 + \sum_{i\in S\setminus\{1\}} (a_i-c_i)r_i,\label{eq3}
\end{eqnarray}
while applying Lemma \ref{lem:distr} to the fourth equation gives
\begin{eqnarray*}  
V_{11} = \frac{1}{\alpha}\Big(a_2b_1c_1r_1' + a_1b_2c_2r_2' +
	\!\! \sum_{i\in S\setminus\{1,2\}} (a_2 b_1 c_i - a_i b_1 c_2 + a_1 c_2 b_i)r_i \Big),
\end{eqnarray*}
where $\alpha=a_2b_1-b_1c_2+a_1c_2$. Substituting according to (\ref{eq2}) and (\ref{eq3}) then yields
\begin{eqnarray*}  
V_{11} < \frac{a_2b_1}{\alpha}\Big(a_1r_1 + \sum_{i\in S\setminus\{1\}} \!\!(a_i-c_i)r_i\Big)
		+ \frac{a_1c_2}{\alpha}\Big(a_2r_2 + \sum_{i\in S\setminus\{2\}} \!\!(a_i-b_i)r_i \Big) \\
	+\frac{a_2 b_1}{\alpha}\!\! \sum_{i\in S\setminus\{1,2\}}\!\!\!\! c_i r_i -
	 \frac{b_1 c_2}{\alpha}\!\! \sum_{i\in S\setminus\{1,2\}}\!\!\!\! a_i r_i   + 
	 \frac{ a_1 c_2}{\alpha}\!\! \sum_{i\in S\setminus\{1,2\}}\!\!\!\! b_i r_i \\
=  \frac{1}{\alpha}\Big( a_1 a_2 b_1 r_1 + a_2 b_1(a_2-c_2) r_2 + a_1 a_2 c_2 r_2+ a_1 c_2 (a_1-b_1) r_1 \\
		+(a_2 b_1 + a_1 c_2 -b_1 c_2)\!\! \sum_{i\in S\setminus\{1,2\}} \!\!\!\! a_i r_i 	 \Big) \\
= \frac{a_2b_1-b_1c_2+a_1c_2}{\alpha}\Big( a_1 r_1 + a_2 r_2 + \!\! \sum_{i\in S\setminus\{1,2\}}\!\!\!\! a_i r_i  \Big) 
= V_{00}.
\end{eqnarray*}
Obviously, replacing `$>$' with `$\geq$', `$<$' or `$\leq$' throughout the proof yields the analogous 
result for the other cases, which finishes the proof.
\end{proof}

The following is a collection of simple consequences of Lemma \ref{lem:comp}.

\begin{corollary}\label{cor:comp}
Let $V_{00},V_{01},V_{10},V_{11}$ and $a,b$ be as in Lemma \ref{lem:comp}. 
Then the following implications hold:
\begin{enumerate}
\renewcommand{\theenumi}{(\roman{enumi})}
\item  \quad$V_{ab} < V_{a,\neg b}, V_{\neg a, b} \quad \Longrightarrow\quad 
							V_{\neg a, \neg b} > \min(V_{a,\neg b},V_{\neg a,b})$.
\item   \quad$V_{ab} > V_{a,\neg b}, V_{\neg a, b}  \quad\Longrightarrow  \quad 
							V_{\neg a, \neg b} < \max(V_{a,\neg b},V_{\neg a,b})$.
\item   \quad$V_{ab} \leq V_{a,\neg b}, V_{\neg a, b}  \quad\Longrightarrow  \quad 
							V_{\neg a, \neg b} \geq \min(V_{a,\neg b},V_{\neg a,b})$.
\item   \quad$V_{ab} \geq V_{a,\neg b}, V_{\neg a, b}  \quad\Longrightarrow  \quad 
							V_{\neg a, \neg b} \leq \max(V_{a,\neg b},V_{\neg a,b})$.
\item    \quad$V_{ab} = V_{a,\neg b} = V_{\neg a, b}  \quad\Longrightarrow  \quad V_{\neg a, \neg b} = V_{ab}$.
\item   \quad $V_{ab}, V_{\neg a, \neg b} \geq  V_{a,\neg b}, V_{\neg a, b}  \quad\Longrightarrow  \quad
V_{ab} = V_{a,\neg b} = V_{\neg a, b}=V_{\neg a, \neg b}$.
\end{enumerate}
\end{corollary}

\begin{proof}
(i), (ii), (iii), (iv) are mere reformulations of Lemma \ref{lem:comp}, while (vi) is an easy consequence.
Thus let us consider (v). If $V_{\neg a, \neg b}$ were $<V_{a,\neg b} = V_{\neg a, b}$, then by Lemma 
\ref{lem:comp} $V_{ab} > \min(V_{\neg a,b}, V_{a,\neg b})$,
contradicting our assumption. Since a similar contradiction crops up if we assume that 
$V_{\neg a, \neg b}>V_{a,\neg b} = V_{\neg a, b}$, it follows that 
$V_{\neg a, \neg b} = V_{a,\neg b} = V_{\neg a, b} =V_{ab}$. 
\end{proof}

Now, in order to prove the theorem, we ignore all states where the optimal policies $\pi_1^\circ,\pi_2^\circ$
coincide. For the remaining $s$ states we denote the actions of $\pi_1^\circ$ by 0 
and those of  $\pi_2^\circ$ by 1. Thus any combination of $\pi_1^\circ,\pi_2^\circ$ can be expressed as a 
sequence of $s$ elements $\in\{0,1\}$, where we assume an arbitrary order on the set of states (take e.g.\ 
the one used in the matrices $P_{ij}$).
We now define sets of policies or sequences, respectively, as follows: First, let $\Theta_i$ be the set of 
policies with exactly $i$ occurrences of 1. Then set $\Pi_0:=\Theta_0=\{00\ldots0\}$, and for $1\leq i\leq s$
\[  \Pi_i:=\{\pi\in\Theta_i \,|\, d(\pi, \pi_{i-1}^*) = 1 \}, \]
where $d$ denotes the Hamming distance, and $\pi_i^*$ is a (fixed) policy in $\Pi_i$ with 
$V(\pi_i^*)=\max_{\pi\in \Pi_i} V(\pi)$. 
Thus, a policy is $\in\Pi_i$, if and only if it can be obtained from $\pi_{i-1}^*$ by replacing a 0 with a 1. \\

\begin{lemma}\label{lem:chain}
$V(\pi_{i-1}^*)\geq V(\pi_i^*)$ for $1\leq i\leq s$. 
\end{lemma}
\begin{proof}
The lemma obviously holds for $i=1$, since $\pi_0^*=00\ldots 0=\pi_1^\circ$ is by assumption optimal. 
Proceeding by induction, let $i>1$ and assume that $V(\pi_{i-2}^*)\geq V(\pi_{i-1}^*)$. 
By construction of the elements in each $\Pi_j$ the policies $\pi_{i-2}^*$, $\pi_{i-1}^*$ 
and $\pi_{i}^*$ differ in at most two states, i.e. the situation is as follows:
\begin{eqnarray*}  
\pi_{i-2}^* &=& \ldots 0 \ldots 0 \ldots\\
\pi_{i-1}^* &=& \ldots 1 \ldots 0 \ldots\\
\pi_{i}^* &=& \ldots 1 \ldots 1 \ldots\\
\pi' &=& \ldots 0 \ldots 1 \ldots
\end{eqnarray*}
Define a policy $\pi'\in\Pi_{i-1}$ as indicated above. Then $V(\pi_{i-2}^*) \geq V(\pi_{i-1}^*) \geq V(\pi')$ 
by induction assumption and optimality of $\pi_{i-1}^*$ in $\Pi_{i-1}$. Applying (iv) of Corollary 
\ref{cor:comp} yields that $V(\pi_i^*)\leq \max(V(\pi_{i-1}^*), V(\pi'))=V(\pi_{i-1}^*)$, which proves 
the lemma.
\end{proof}

Since the policies $\pi_0^*=00\ldots 0=\pi_1^\circ$ and $\pi_s^*=11\ldots 1=\pi_2^\circ$ are assumed to be
optimal, it follows that all policies $\pi_i^*$ are optimal as well.  Now we are able to prove the Theorem 
by  induction on the number of states $s$ where the policies $\pi_1^\circ, \pi_2^\circ$ 
differ. For $s=1$ it is trivial, while for $s=2$ there are two combinations of $\pi_1^\circ=\pi_0^*$ and 
$\pi_2^\circ=\pi_2^*$. One of them is identical to $\pi_1^*$ and hence optimal, while the other one is 
optimal due to Corollary \ref{cor:comp} (v).

Thus, let us assume that $s>2$. Then we have already shown that the policies $\pi_i^*$ and hence in particular
$\pi_1^*=00\ldots 010\ldots 0$ and $\pi_{s-1}^*=11\ldots 101\ldots 1$ are optimal. Since $\pi_1^*$ and
$\pi_s^*=11\ldots 1=\pi_2^\circ$ are optimal policies that share a common digit in position $k$, we may conclude 
by induction assumption that all policies with a 1 in position $k$ are optimal. A similar argument applied to the
policies  $\pi_{0}^*=00\ldots 0$ and $\pi_{s-1}^*$ shows that all policies with a 0 in position $\ell$ (the position of the 0
in $\pi_{s-1}^*$) are optimal. Note that by construction of the sets $\Pi_i$, $k\neq \ell$. Thus, we have shown
that all considered policies are optimal, except those with a 1 in position $\ell$ and a 0 in position $k$. 
However, as all policies of the form
\begin{eqnarray*}  
&& \,\,\quad k \quad \,\,\,\ell\\
&& \ldots 0 \ldots 0 \ldots\\
&& \ldots 1 \ldots 0 \ldots\\
 && \ldots 1 \ldots 1 \ldots
\end{eqnarray*}
are optimal, a final application of Corollary \ref{cor:comp} (v) shows that these are optimal as well.\qed

\section{Mixing Optimal Policies}
Theorem \ref{theorem} can be extended to \textit{mixing} optimal policies, that is, our policies are 
not deterministic (\textit{pure}) anymore, but in each state we choose an action randomly. Building 
up on Theorem \ref{theorem} we can show that any mixture of optimal policies is optimal as well.

\begin{theorem}\label{thm2}
Let $\Pi^*$ be a set of pure optimal policies on a unichain MDP $\mathcal M$. Then any policy that chooses
at each visit in each state $i$ randomly an action $a$ such that there is a policy $\pi\in \Pi^*$ with $a=\pi(i)$, 
is optimal.
\end{theorem}

The theorem will be obtained with the help of the following Lemma.
\begin{lemma}\label{lem:mix}
Let $\pi_1,\pi_2$ be two policies on a unichain MDP $\mathcal M$ that differ only in a single state $s_1$,
i.e.\  $\pi_1(i)=\pi_2(i)$ for all $i\neq s_1$ and $\pi_1(s_1)\neq\pi_2(s_1)$. 
Let $\mu_1=(a_i)_{i\in S}$, $\mu_2=(b_i)_{i\in S}$ be the invariant distributions and $V_1,V_2$ the average 
rewards of $\pi_1$ and $\pi_2$, respectively. Then the mixed policy $\pi$ that chooses in $s_1$ action
$\pi_1(s_1)$ with probability $\lambda$ and $\pi_2(s_1)$ with probability $(1-\lambda)$ and coincides
with $\pi_i$ in all other states has invariant distribution $\mu=(c_i)_{i\in S}$ with
\[  c_i = \frac{\lambda a_i b_{s_1} + (1-\lambda)a_{s_1}b_i}{\lambda b_{s_1}+(1-\lambda)a_{s_1}} \]
and average reward 
\[  V=\frac{\lambda b_1 V_1 + (1-\lambda)a_1V_2}{\lambda b_1+(1-\lambda)a_1}.  \]
\end{lemma}

\begin{proof}
First, note that the transition matrices $P_1,P_2$ of $\pi_1,\pi_2$ and $P$ of $\pi$ share all
rows except row $s_1$, which we assume to be the first row. Furthermore we write $p_{ij}$ for
$p_{\pi_1(i)}(i,j)$ and $q_{1j}$ for $p_{\pi_2(s_1)}(s_1,i)$, so that the entries of row $1$ in $P$ 
are of the form $\lambda p_{1j}+ (1-\lambda)q_{1j}$.
Now, let $\nu:=(\nu_i)_{i\in S}$ with $\nu_i=\lambda a_i b_1 + (1-\lambda)a_1b_i$. Then
\begin{eqnarray*}  
(\nu P)_i &=& a_1 b_1(\lambda p_{1i} + (1-\lambda)q_{1i}) +  
		\sum_{j\in S\setminus\{1\}}  \big( \lambda a_j b_1 +  (1-\lambda)a_1 b_j \big) p_{ji}=\\
&=& \lambda b_1  \big(a_1 p_{1i} + \sum_{j\in S\setminus\{1\}} a_j p_{ji} \big) + 
		(1-\lambda) a_1  \big(b_1 q_{1i} + \sum_{j\in S\setminus\{1\}} b_j  p_{ji}\big)  = \\
&=& \lambda a_i b_1 + (1-\lambda)a_1b_i = \nu_i,
\end{eqnarray*}  
since the $a_i$'s and $b_i$'s form an invariant distribution of  $P_1,P_2$, respectively. Since the 
$c_i$'s are only a normalized version of the $\nu_i$'s, this finishes the first part of the proof.

Now, given the invariant distribution of $\pi$, its average reward can be written as $V= \sum_{i\in S} c_i r_{\pi(i)}$.
Thus, writing $r_i$ for $r_{\pi_1(i)}$ and $r_1'$ for $r_{\pi_2(s_1)}$ one has
\begin{eqnarray*}  
V &=& \frac{a_1 b_1}{\lambda b_1 +(1-\lambda) a_1}\big( \lambda r_1 + (1-\lambda)r_1' \big)  +
	\sum_{i\in S\setminus\{1\}} \frac{\lambda a_i b_1 + (1-\lambda)a_1b_i}{\lambda b_1+(1-\lambda)a_1}  r_i=\\
&=&  \frac{1}{\lambda b_1 +(1-\lambda) a_1} 
		\Big( \lambda b_1 \sum_{i\in S} a_i r_i +  
			(1-\lambda) a_1 \big( b_1 r_1'  +  \sum_{i\in S\setminus\{1\}} b_i r_i \big)  \Big) =\\
&=& \frac{\lambda b_1 V_1 + (1-\lambda)a_1V_2}{\lambda b_1+(1-\lambda)a_1}. 
\end{eqnarray*}  
\end{proof}

\begin{proof}[Proof of Theorem \ref{thm2}]
Let us first assume the simplest case, where we have two optimal policies $\pi_1^\circ,\pi_2^\circ$ that differ 
only in some state single $s_1$. By Lemma \ref{lem:mix}, the policy resulting from $\pi_1^\circ$ and $\pi_2^\circ$  
when mixing actions in state $s_1$ has the same average reward as $\pi_1^\circ$ and $\pi_2^\circ$ and therefore
is optimal. Now, each mixture of actions in a single state can be interpreted as a new action in this state.
Thus, proceeding by induction, the mixture of $n$ optimal policies $\pi_1^\circ,\ldots,\pi_n^\circ$ that 
differ only in a single state is optimal as well.

Now, in the general case, where we want to mix actions in $s>1$ states, we have at each state
$i$ the actions (of some pure optimal policies $\in\Pi^*$) $a_1^i,a_2^i,\ldots,a_{k_i}^i$ at our disposal.
By Theorem \ref{theorem} all combinations $(a_{j_1}^1, a_{j_2}^2,\ldots,a_{j_s}^{s})$ with $1\leq j_i\leq k_i$ 
for all $i$ are optimal as well. Thus, we may fix the actions in $s-1$ states so that we have e.g.\ optimal policies
of the form $(a_{j}^1, a_{1}^2,\ldots,a_{1}^{s})$ with $1\leq j \leq k_1$. As we have seen above, 
all policies that are obtained by mixing all available actions in the first state are optimal. Furthermore, each 
mixture can again be interpreted as new available action, so that we may repeat our argument for each of
the remaining states, thus showing that each mixed optimal policy is optimal, too.

So far, we have considered only the case where the relative frequencies with which the actions in a fixed
state are chosen converge. If this does not hold, it may happen that the process does not converge to an 
invariant distribution. However, the average rewards after $t$ steps converge nevertheless. 
Let $\lambda_i^t(a)$ be the relative frequency with which action $a$ was chosen in state $i$ after
$t$ steps in $i$, and let $\mu_t$ be the distribution over the states after these $t$ steps. 
Then the average reward $V_t$ thereby obtained is $\sum_{i\in S} \mu_t(i) \sum_a \lambda_i^t(a)  r_a(i)$. 
This is of course also the expected average reward after $t$ steps when constantly choosing action $a$ in
state $i$ with probability $\lambda_i(a):=\lambda_i^t(a)$ for each $i,a$. As each of these sequences
has already been shown to converge to the optimal value $V^*$, we have the following situation.
For each $V_{t_1}$ of the sequence $(V_t)_{t\in\mathbb N}$ there is a sequence $(V_t(\pi))_{t\in\mathbb N}$ 
with $\lim_{t\to \infty} V_t(\pi)=V^*$ such that $V_{t_1}(\pi)=V_{t_1}$. It follows that 
$\lim_{t\to \infty} V_t=V^*$.
\end{proof}

\section{Extensions, Applications and Remarks}
\subsection{Optimality is Necessary}
Given some policies with equal average reward $V$, in general, it is not the case 
that a combination of these policies again has average reward $V$, as the following example shows.
Thus, optimality is a necessary condition in Theorem \ref{theorem}.
\begin{exa}\label{ex2}
Let $S=\{s_1,s_2\}$ and $A=\{a_1,a_2\}$. The transition probabilities are given by
\begin{eqnarray*}  
(p_{a_1}(i,j))_{i,j\in S}=(p_{a_2}(i,j))_{i,j\in S}=
\left(\begin{array}{*{2}{cc}}
0 & 1 \\
1 & 0 \\
\end{array}\right),
\end{eqnarray*}  
while the rewards are $r_{a_1}(s_1)=r_{a_1}(s_2)=0$ and $r_{a_2}(s_1)=r_{a_2}(s_2)=1$. 
Since the transition probabilities of all policies are identical, policy $(a_2,a_2)$ with an 
average reward of 1 is obviously optimal. Policy $(a_2,a_2)$ can be obtained as a combination
of the policies $(a_1,a_2)$, and $(a_2,a_1)$, which however only yield an average reward of 
$\frac{1}{2}$.
\end{exa}

\subsection{Multichain and Infinite MDPs}
Theorem \ref{theorem} does not hold for MDPs that are not unichain as the following simple example demonstrates.
\begin{exa}\label{ex}
Let $S=\{s_1,s_2\}$ and $A=\{a_1,a_2\}$. The transition probabilities are given by
\begin{eqnarray*}  
(p_{a_1}(i,j))_{i,j\in S}=
\left(\begin{array}{*{2}{cc}}
1 & 0 \\
0 & 1 \\
\end{array}\right),\quad
(p_{a_2}(i,j))_{i,j\in S}=
\left(\begin{array}{*{2}{cc}}
\frac{1}{2} & \frac{1}{2} \\
\frac{1}{2} & \frac{1}{2} 
\end{array}\right),
\end{eqnarray*}  
while the rewards are $r_{a_1}(s_1)=r_{a_1}(s_2)=1$ and $r_{a_2}(s_1)=r_{a_2}(s_2)=0$. Then
the policies $(a_1,a_1), (a_1,a_2), (a_2,a_1)$ all gain an average reward of 1 and are optimal,
while the combined policy yields suboptimal average reward 0.
\end{exa}

Even though this seems to be quite a strict counterexample (note that the MDP is even communicating), 
we think that in certain restricted settings Theorems \ref{theorem} and \ref{thm2} will hold as well. For example,
adding a set of states that are transient under every policy does not matter. Furthermore, if the components
of a multichain MDP are the same under every policy, it is obvious that the Theorems hold as well. However,
things become more complicated, if the set of transient states or the components change with the policy
as in Example \ref{ex}. Nevertheless, extensions of our results to the multichain case don't seem to be 
impossible as such, but may work under some clever restrictions, e.g. by combining exclusively in states 
that are not transient under any policy. In any case the main task when working on such extensions will
probably be to determine what policy changes will result in what changes in the set of transient states and
components, respectively.\\

The situation for MDPs with countable set of states/actions is similar. Under the (strong) 
assumption that there exists a unique invariant and positive distribution for each policy,
Theorems \ref{theorem} and \ref{thm2} also hold for these MDPs. In this case the proofs are 
identical to the case of finite MDPs (with the only difference that the induction becomes transfinite).
However, in general, countable MDPs are much harder to handle as optimal policies need not
be stationary anymore (cf. \cite{pute}, p.413f).

\subsection{An Application}
Even though the presented results may seem more of theoretical interest, there is a straightforward
application of Theorem \ref{thm2}, which actually was the starting point of this paper. Consider
an algorithm operating on an MDP that every now and then recalculates the optimal policy 
according to its estimates of the transition probabilities and the rewards, respectively. Sooner
or later the estimates are good enough so that the calculated policy is indeed an optimal one.
However, if there is more than one optimal policy, it may happen that the algorithm does not 
stick to a single optimal policy but starts mixing optimal policies irregularly. Theorem \ref{thm2} 
guarantees that the average reward of such a process again is still optimal.

\section{Conclusion}
We conclude with a more philosophical remark. MDPs are usually presented as a standard example
for decision processes with delayed feedback. That is, an optimal policy often has to accept locally 
small rewards in present states in order to gain large rewards later in future states. One may think 
that this induces some sort of context in which actions are optimal, e.g. that choosing a locally suboptimal
action only ``makes sense'' in the context of heading to the higher reward states. Our results however show
that this is not the case and optimal actions are rather optimal in any context.


\end{document}